\documentclass[11pt]{article}

\usepackage[a4paper,margin=1in]{geometry}
\usepackage{comment}
\usepackage{amsmath, amssymb, amsthm, mathtools}
\usepackage{mathrsfs}
\usepackage{enumitem}
\usepackage{hyperref}
\usepackage{cleveref}
\usepackage{tikz}
\usepackage{tikz-cd}
\usepackage{authblk}

\newtheorem{theorem}{Theorem}[section]
\newtheorem{proposition}[theorem]{Proposition}
\newtheorem{lemma}[theorem]{Lemma}
\newtheorem{corollary}[theorem]{Corollary}

\theoremstyle{definition}
\newtheorem{definition}[theorem]{Definition}

\newtheorem{fact}[theorem]{Fact}

\newcommand{\Homeo}{\operatorname{Homeo}}

\newcommand{\Ends}{\operatorname{Ends}}

\DeclareMathOperator{\Map}{Map}
\DeclareMathOperator{\Stab}{Stab}
\newcommand{\genus}{\operatorname{genus}}
\newcommand\restr[2]{{
  \left.\kern-\nulldelimiterspace 
  #1 
  \littletaller 
  \right|_{#2} 
  }}

\newcommand{\littletaller}{\mathchoice{\vphantom{\big|}}{}{}{}}
 
\makeatletter
\def\@maketitle{%
  \newpage
  \null
  \vskip 2em%
  \begin{center}%
  \let \footnote \thanks
    {\LARGE \bfseries \@title \par}
    \vskip 1.5em%
    {\large
      \lineskip .5em%
      \begin{tabular}[t]{c}%
        \@author
      \end{tabular}\par}%
    \vskip 1em%
    {\large \@date}%
  \end{center}%
  \par
}
\makeatother

\title{Finitely supported mapping classes
on infinite-type surfaces}
\author[]{Ekaterina Kuritsina}
\date{}

\begin{document}
\maketitle
\begin{quote}
\small 
\noindent We study the topological closure of the group of mapping classes group with finite-type support on an infinite-type surface. Since every such mapping class fixes each end which is not a puncture, we have that the closure is contained in the pointwise stabilizer of infinite-type ends. A natural question is therefore when this inclusion is an equality, and if it is not, what must be added to obtain the full stabilizer.  We prove that finitely supported mapping classes are dense in the stabilizer precisely when there is at most one end accumulated by genus and at most one end accumulated by punctures. We show that if we remove either of these conditions, we need to add puncture and/or handle shifts to topologically generate the stabilizer.
\end{quote}


\section{Introduction}
Let $S$ be an orientable connected surface, possibly with boundary.
We will consider its mapping class group $\Map(S)=\Map(S,\partial S):= \pi_0(\Homeo^{+}(S,\partial S))$, the group of orientation-preserving homeomorphisms up to homotopy, fixing the boundary point-wise. It is endowed with the topology induced by the compact-open topology on $ \Homeo^{+}(S,\partial S)$. 

For mapping class groups of finite-type surfaces, generating sets are well understood. In particular, Dehn twists about suitable collections of curves provide finite generating sets for mapping class groups of closed, orientable surfaces \cite{Lickorish1964}. More generally, finite generating sets are known for mapping class groups of finite-type surfaces; see, for example, Farb and Margalit \cite{FarbMargalit}.

The situation is more difficult for surfaces of infinite type. Their mapping class groups are uncountable, and the subgroup of mapping classes with finite-type support, denoted by $\Map_f(S)$, is a proper subgroup. Still, it might be dense, and our main objective is to understand when that is the case.

We determine when $\Map_f(S)$ coincides with the point-wise stabilizer of the space of ends after removing punctures, and identify the shifts that have to be added to obtain the full stabilizer in the remaining cases. 
More precisely, set
\[
 \qquad \tilde{E}:=\Ends(S)\setminus\{\text{punctures}\}, \ \operatorname{Stab}(\tilde{E})=\bigcap_{e\in \widetilde{E}}\operatorname{Stab}(e).
\]

Our main result is the following.

\begin{theorem}
\label{maintheorem}

Let $S$ be a surface with empty boundary. The point-wise stabilizer of $\tilde{E}$ equals to the topological closure of $\Map_f(S)$
if and only if there exists at most one end accumulated by punctures and at most one end accumulated by genus.\\
\\
Furthermore:
\begin{itemize}
    \item []
    If there exists at most one end accumulated by punctures, then
    $$
    \Stab(\tilde{E})
    =
    \overline{\langle \Map_f(S),\{\text{handle shifts}\}\rangle}.
    $$

    \item []
    If there exists at most one end accumulated by genus, then
     $$
    \Stab(\tilde{E})
    =
    \overline{\langle \Map_f(S),\{\text{puncture shifts}\}\rangle}.
 $$

    \item [] For every surface,
     $$
    \Stab(\tilde{E})
    =
    \overline{\langle
    \Map_f(S),
    \{\text{handle shifts}\},
    \{\text{puncture shifts}\}
    \rangle}.
     $$
\end{itemize}
\end{theorem}

Shifts --- homeomorphisms that translate a periodic pattern (handles, punctures, or more general subsurfaces) along a surface --- appear repeatedly throughout the study of the big mapping class groups, which is part of why they are a natural object to understand on their own. 
In particular Abbott, Miller, and Patel \cite{AbbottMillerPatel} used compositions of shift maps --- together with a pseudo-Anosov on a finite-type subsurface --- to construct the first known loxodromic isometries of the relative arc graph for a broad class of infinite-type surfaces, extending an earlier construction of Bavard \cite{Bavard} for the plane minus a Cantor set. In follow-up work \cite{AbbottMillerPatelShifts}, the same authors show that composing a Dehn twist with a shift map already produces a loxodromic isometry. Lanier and Loving \cite{LanierLoving} used shifts to construct subgroups isomorphic to the restricted wreath product $\mathbb{Z}\wr\mathbb{Z}$. Their construction shows that every mapping class group of infinite-type surface contains a subgroup that is not virtually abelian and contains no nonabelian free subgroup, and hence does not satisfy the strong Tits alternative. They also exhibited examples that fail the classical Tits alternative.

In \cite{PatelVlamis}, Patel and Vlamis considered a similar question to the one in our paper. They studied the pure mapping class group $\operatorname{PMap}(S)$ and determined a topological generating set in terms of Dehn twists and handle shifts. In particular, when $S$ has at most one end accumulated by genus, Dehn twists topologically generate $\operatorname{PMap}(S)$. When $S$ has more than one end accumulated by genus, one needs to add handle shifts to the set of Dehn twists. 
The proof of Theorem \ref{maintheorem} is inspired by the work of Patel and Vlamis.
Their result was improved by Aramayona, Patel, and Vlamis \cite{AramayonaPatelVlamis}.
They prove that
\[
\operatorname{PMap}(S) = \overline{\operatorname{PMap}_c(S)} \rtimes \kappa\big(H^1_{\mathrm{sep}}(\hat S,\mathbb{Z})\big),
\]
where $\operatorname{PMap}_c(S)$ is the compactly supported mapping class group and $\kappa\big(H^1_{\mathrm{sep}}(\hat S,\mathbb{Z})\big)$ is a free abelian subgroup generated by handle shifts. This group has one generator for each end of $S$ accumulated by genus, when the number of the ends is finite and greater than one; if $S$ has infinitely many such ends, it is isomorphic to the Baer--Specker group $\mathbb{Z}^{\mathbb{N}}$, the group of all integer sequences with coordinate-wise addition. Note that other authors considered problems related to algebraic and normal generation in the mapping class group of infinite-type surfaces (see \cite{DomatFanoniHensel} and \cite{MalesteinTao}).

It is natural to try to compare our result with that of Patel and Vlamis \cite{PatelVlamis} and understand if and when we can replace $\Map_f(S)$ with $\Map_c(S)$ in our theorem. Clearly, $\overline{
    \Map_f(S)} \subseteq \overline{
    \Map_c(S)} $, and they coincide if $S$ has no punctures. In general, we show:
\begin{lemma}
\label{prop}
 The surface $S$ has at least one end accumulated by punctures if and only if
\[
\overline{
\left\langle
\operatorname{Map}_f(S),\{\text{puncture shifts}\}
\right\rangle}
=
\overline{
\left\langle
\operatorname{Map}_c(S),\{\text{puncture shifts}\}
\right\rangle}.
\]
\end{lemma}
As a consequence, we can restate part of our result in terms of $\Map_c(S)$.
\begin{corollary}
\label{cor1}
If there is at least one end accumulated by punctures, then:
\begin{itemize}

    \item []
    If there exists at most one end accumulated by genus, then
     $$
    \Stab(\tilde{E})
    =
    \overline{\langle \Map_c(S),\{\text{puncture shifts}\}\rangle}.
 $$

    \item [] In general,
     $$
    \Stab(\tilde{E})
    =
    \overline{\langle
    \Map_c(S),
    \{\text{handle shifts}\},
    \{\text{puncture shifts}\}
    \rangle}.
     $$
\end{itemize}
\end{corollary} 

A natural application of our results concerns the infinite braid group. The relationship between the group and the mapping class group has also been considered by other authors: \cite{LeRoux1997}, \cite{Wiest2003}, \cite{Fabel2003}, \cite{Brunner2008}. Recall that the classical braid group \(B_n\) can be identified with the mapping class group of a disk with \(n\) punctures, denoted by $\mathbb{D}_n$. The infinite braid group is defined as the direct limit of the classical ones, with respect to the standard inclusion $B_n \hookrightarrow B_{n+1}$:
$$
B_\infty=\varinjlim B_n=\varinjlim \operatorname{Map}(\mathbb{D}_n)
$$

Let $
S=\mathbb{C}\setminus\mathbb{N}
$
and, for every \(n\geq 1\), let
$
D_n=\overline{B}\left(\frac{n+1}{2},\frac{n+1}{2}\right)
$ be a closed disk with center at $\frac{n+1}{2}$ and radius $\frac{n+1}{2}$.
We identify $\mathbb{D}_n$ with \(D_n\cap S\), so that $B_\infty$ is identified with \(\Map_f(S)\). Our main theorem implies that the density of \(B_\infty\) in \(\operatorname{Map}(S)\).
\begin{corollary}
With the notation above,
$$
\overline{B_\infty}
=
\operatorname{Map}(S)
.
$$
\end{corollary}
The paper is organized as follows.
In \Cref{sec:preliminaries}, we recall the necessary background
on mapping class group and known results which will be used to prove the theorem.
In \Cref{sec:main-result}, we prove our main result.

\textbf{Acknowledgments.} I am grateful to Federica Fanoni, the author's PhD advisor, for her support and guidance throughout this project. In particular, I thank her for suggesting the project, providing valuable guidance at every stage, and helping me identify and correct  mistakes. The author was supported by the ANR PRC grant GALS (reference ANR-23-CE40-0001).

\section{Preliminaries}
\label{sec:preliminaries}

Throughout the paper, $S$ denotes an orientable connected surface, possibly with boundary. We fix the notation $S\setminus \alpha $ to mean cutting \(S\) along \(\alpha\), so that the resulting surface has two boundary components corresponding to the two sides of $\alpha$.
\subsection{Space of Ends}
\label{subsec:ends}
We will just state some definitions here. For more details about properties of end spaces and mapping class groups of infinite-type surfaces in general, we refer the reader to Aramayona's and Vlamis' paper \cite{Aramayona2020}.
Let $S$ be a connected surface. Choose an exhaustion
\[
    K_1 \subset K_2 \subset \cdots \subset S
\]
by compact subsurfaces.
For each $n$, let $\mathcal{C}_n$ denote the set of connected components of
$S\setminus K_n$. Since $K_n\subset K_{n+1}$, every component of
$S\setminus K_{n+1}$ is contained in a unique component of $S\setminus K_n$.
Thus, for each $n$ there is a natural inclusion
\[
    \pi_n:\mathcal{C}_{n+1}\longrightarrow\mathcal{C}_n.
\]

\begin{definition}
The \emph{space of ends} of $S$, denoted by $\Ends(S)$, is the inverse limit
\[
    \Ends(S)
    :=
    \varprojlim
    \left(
        \mathcal{C}_n,\pi_n
    \right),
\]
with the induced topology.
Thus, an end of $S$ is a sequence
\[
    e=(C_1,C_2,\ldots)
\]
such that $C_n\in\mathcal{C}_n$ and
\[
    C_{n+1}\subset C_n
\]
for every $n$.
\end{definition}

If $U\subset S$ is an open subset,
we set
\[
    U^*
    :=
    \left\{
        e=(C_1,C_2,\ldots)\in\Ends(S)
        \,\middle|\,
        C_n\subset U
        \text{ for some sufficiently large }n
    \right\}.
\]
It can be shown that the collection of sets of the form $U^*$, for $U$ varying among the open subsets of $S$, forms a basis for the topology of $\Ends(S)$. 

Throughout the paper, planar isolated ends will be referred to as \emph{punctures}.
\begin{definition}
Let $e\in\Ends(S)$. An open subset $U\subset S$ with compact boundary is
called a \emph{neighborhood of $e$} if
\[
    e\in U^*.
\]
\end{definition}
\begin{definition}
An end is called \emph{planar} if it has a neighborhood of genus zero, and \emph{non-planar} or \emph{accumulated by genus} otherwise. 
\end{definition}
\begin{definition}
An end $e\in\Ends(S)$ is said to be \emph{accumulated by punctures} if every neighborhood of $e$ contains infinitely many punctures.
\end{definition}
We set
$$
    \widetilde E
    :=
    \Ends(S)\setminus\{\text{punctures}\}.
$$

We will be interested in the pointwise stabilizer

$$
    \Stab(\widetilde E)
    :=
    \bigcap_{e\in\widetilde E}\Stab(e),
$$
where $\Stab(e)$ denotes the stabilizer of the end $e$ under the natural action of $\Map(S)$ on $\Ends(S)$.

\subsection{Mapping class groups}
We endow $\Map(S)$ with the quotient topology induced by the compact-open topology on $\Homeo^+(S)$. This topology agrees with the permutation topology on $\Map(S)$. A basis for the latter is given by the sets

$$
    \theta U_A
    =
    \left\{
        \theta\psi \,\middle|\,
        \psi(a)=a \text{ for every } a\in A
    \right\},
$$
where $\theta\in\Map(S)$ and $A$ is a finite collection of homotopy classes of essential simple closed curves on $S$.
\begin{definition}
We say that a subset $X\subseteq G$ \emph{topologically generates} a topological group $G$ if

$$
\overline{\langle X\rangle}=G,
$$
where $\langle X\rangle$ denotes the subgroup of $G$ generated by $X$.
\end{definition}
\begin{definition}
A mapping class $\phi \in\Map(S)$ is called \emph{compactly supported} if it has a representative with compact support. We denote the subgroup of compactly supported mapping classes by $\Map_c(S).$

A mapping class $\phi \in\Map(S)$ is called \emph{finitely supported} if it has a representative supported on a finite-type subsurface of $S$. We denote the subgroup of finitely supported mapping classes by $\Map_f(S).$
\end{definition}
Clearly,

$$
    \Map_c(S)\leq \Map_f(S).
$$

We will use the following standard consequence of the classification of surfaces of finite type.

\begin{fact}
\label{classification}
Let $S$ be a finite type surface, possibly with boundary and let $\alpha, \beta \subset S$ be homotopy classes of essential simple closed
curves. The following are equivalent:
\begin{enumerate}
  \item The homotopy classes $\alpha$ and $\beta$ belong to the same $\mathrm{Map}(S)$-orbit.
  \item Either both representative curves are non-separating or there is a bijection $\varphi$ between the set of connected components of
  $S \setminus \alpha$ and the set of connected components of $S \setminus \beta$ such that
  for every component $W$ of $S\setminus \alpha $ ,
  $
    \mathrm{genus}(W) = \mathrm{genus}(\varphi(W)), \
    \partial W = \partial \varphi(W),
  $
  and $\# \text{punctures}(W)=\# \text{punctures}(\varphi(W))$.
\end{enumerate}
\end{fact}
A \emph{multicurve} on a surface $S$ is a finite collection of pairwise disjoint, pairwise non-isotopic, essential simple closed curves contained in $S$.
We also record the following for later use:
\begin{fact}
\label{lem:end-splitting}
Let $\Gamma\subset S$ be a multicurve, and let $
\{W_j\}_{j\in J}
$
be the connected components of the surface obtained by cutting $S$ along $\Gamma$. Then every end of $W_j$ determines an end of $S$, giving an injective map

$$
\Ends(W_j)\longrightarrow\Ends(S).
$$
Moreover, the images of these maps for distinct components $W_j$ are pairwise disjoint.
\end{fact}
The following observation will be useful when constructing finite-type subsurfaces whose complementary components contain prescribed sets of ends.

\begin{fact}
\label{fact:cantor-separating-curve}
Let $S$ be any surface with $\Ends(S)= U \sqcup V$, where $U, V$ are non-empty, clopen subsets. Then there exists an essential simple closed curve
$\gamma\subset S$ such that  $S\setminus\gamma=\Sigma_1 \cup\Sigma_2$ and $\Ends({\Sigma_1})= U, \ \Ends(\Sigma_2)= V$.
\end{fact}

We will also use the standard cutting homomorphism for mapping class groups of finite-type surfaces, see for the reference \cite{FarbMargalit}.

\begin{fact}[Cutting homomorphism]
\label{fact:cutting-homomorphism}
Let $S$ be a closed surface with finitely many marked points, and let
$a_1,\ldots,a_n$ be pairwise disjoint, pairwise non-isotopic, essential simple closed curves. Let

$$
    \Map(S,\{[a_1],\ldots,[a_n]\})
$$
denote the subgroup of mapping classes preserving the set of isotopy classes
${[a_1],\ldots,[a_n]}$. Then there is a well-defined homomorphism

$$
    \zeta:
    \Map(S,\{[a_1],\ldots,[a_n]\})
    \longrightarrow
    \Map\left(S\setminus\bigcup_{i=1}^n a_i\right),
$$
whose kernel is generated by the Dehn twists about the curves $a_1,\ldots,a_n$:

$$
    \ker(\zeta)
    =
    \langle T_{a_1},\ldots,T_{a_n}\rangle.
$$

\end{fact}
\subsection{Shifts}
\label{subsec:shifts}
We finally recall the notion of a shift. Let

$$
    \mathbb R\times[-1,1]
$$
be a strip. Remove, for each $n\in\mathbb Z$, an open disk of radius $1/4$ centered at $(n,0)$, and attach a fixed connected surface $\Pi$ with a single boundary component, which is compact, to each resulting boundary circle. We denote the resulting surface by $D_\Pi$.
\begin{definition}
A \emph{shift} on $D_\Pi$ is a homeomorphism which agrees with the translation

$$
    (x,y)\longmapsto(x+1,y)
$$
on a smaller substrip

$$
    (\mathbb R\times[-1+\varepsilon,1-\varepsilon]) \cap D_\Pi
$$
for some $\varepsilon>0$, and which is the identity on the boundary of the strip.
\end{definition}
Suppose that $D_\Pi$ admits a proper embedding into $S$. The shift on $D_\Pi$, extended by the identity outside the image of $D_\Pi$, determines a mapping class of $S$, which we also call a \emph{shift}.

Two special cases will be of particular importance. If $\Pi$ is a once-punctured disk, we call the resulting mapping class a \emph{puncture shift}. If $\Pi$ is a one-holed torus, we call it a \emph{handle shift}.


\section{Main Result}
\label{sec:main-result}

In this section, we prove Theorem \ref{maintheorem}. We first show that the closure of the finitely-supported mapping class group is contained in the stabilizer of all infinite-type ends.
\begin{lemma}
$\overline{\Map_f(S)} \subseteq \operatorname{Stab}(\tilde{E})$.
\end{lemma}\begin{proof}
First, we show that the finitely supported mapping class group is contained in $\operatorname{Stab}(\tilde{E})$. Indeed, a finitely supported mapping class can only move ends contained in its support, which can only be punctures. Therefore, every end of $\widetilde{E}$ is fixed, and hence every finitely supported mapping class belongs to $\operatorname{Stab}(\tilde{E})$.

Next, we show that \(\operatorname{Stab}(\widetilde{E})\) is closed. For each end \(e\), the condition that a mapping class fixes \(e\) defines a closed subset of the mapping class group. Since \(\operatorname{Stab}(\widetilde{E})\) is an intersection of such closed subsets, it is itself closed.
We conclude that
\[
\overline{\operatorname{Map}_f(S)}\subseteq \operatorname{Stab}(\tilde{E}),
\]
as claimed.
\end{proof}
We denote by \textnormal{(P)} and \textnormal{(G)} the following conditions:
\begin{align*}
\textnormal{(P)}\quad &\text{There is at most one end accumulated by punctures},\\
\textnormal{(G)}\quad &\text{There is at most one end accumulated by genus}.
\end{align*}

The next step is to consider a mapping class in the stabilizer and determine under what conditions it can be expressed as a limit of finitely supported mapping classes.

Take $[f] \in \operatorname{Stab}(\tilde{E})$.
Fix a pants decomposition $\{P_1, P_2,...\}$ of S, such that $R_n=\bigcup_{i=1}^{n}P_i$ is connected. Choose a component $a$ of $\partial{P_1}$.
\begin{lemma} There always exists an essential finite-type subsurface $F$ such that:
\begin{itemize}
    \item The curves $a$ and $f(a)$ are contained in $F$,
    \item Each component of $\partial F$ is separating,
    \item There is no component of $S \setminus F$ of finite type.

\end{itemize} If both conditions $(P)$ and $(G)$ hold, then F can additionally be chosen to satisfy the properties:
    \begin{enumerate}[label=(\alph*)]
        \item at most one component of $S \setminus F$ is of positive genus;
        \item at most one component of $S \setminus F$ has punctures.
    \end{enumerate}
If condition $(P)$ holds, then F can be chosen to have property (b). If condition $(G)$ holds, F can be chosen to have property (a). 
\end{lemma}
\begin{proof}
We can always find an essential finite-type subsurface $F_0$ containing $a$ and $f(a)$,
with each component of $\partial F_0$ separating. Since $\partial F_0$ is finite, $S \setminus F_0$ has only finitely many connected
components $W_1, \dots, W_k$. \\
Suppose both conditions hold. Then at most one $W_i$ has infinite genus. We call
this component $W_{i_0}$, if it exists. Likewise, at most one $W_i$ has infinitely
many punctures; we  call it $W_{j_0}$, if it exists.
We now construct $F$ from $F_0$ in two steps.

\emph{Step 1.} Glue into $F_0$ every $W_i$ that is of finite type.

\emph{Step 2.} For every remaining component $W_i$, with $i \neq i_0, j_0$ --
necessarily of finite genus $g_i$ and finitely many punctures $p_i$ -- choose a separating curve in $W_i$ isolating a
finite-type subsurface $W_i'$ containing all the genus and all the punctures
of $W_i$. Such a curve exists because of Fact \ref{fact:cantor-separating-curve}. The complementary piece $W_i \setminus W_i'$ then has genus $0$
and no punctures. We glue $W_i'$ into $F_0$, obtaining a new surface $F$, which satisfies our requirements.\\
If instead only $(P)$ holds, the same argument (omitting the genus step)
shows $F$ satisfies (b); if only $(G)$ holds, the symmetric argument shows
$F$ satisfies (a).
\end{proof}

\begin{lemma} Suppose \(a\) is separating. Then \(a\) and \(f(a)\) induce the same partition of \(\partial F\).\end{lemma}
\begin{proof}
We first prove that the construction of \(F\) guarantees that \(a\) and \(f(a)\) induce the same partition of the boundary components of \(F\).
For every boundary component \(b\) of \(\partial F\), let \(C_b\) be the complementary component of \(S\setminus F\) adjacent to \(b\). Since \(C_b\) is of infinite type, we may choose
$$
e_b\in \Ends(C_b)\cap\widetilde E.
$$
Since \(f\in\operatorname{Stab}(\widetilde E)\), we have
$$
f(e_b)=e_b.
$$
Denote by $U_+$ and $U_-$ connected subsurfaces of $S$ such that

$$
S\setminus a=U_+\sqcup U_-
$$

and

$$
S\setminus f(a)=f(U_+)\sqcup f(U_-).
$$
By Fact ~\ref{lem:end-splitting}, applied to \(a\), the component \(C_b\) is entirely contained in one of \(U_+\) or \(U_-\). Suppose, without loss of generality, that
$
C_b\subset U_+.
$
Then
$$
\Ends(C_b)\subset\Ends(U_+).
$$
Since \(f\) is a homeomorphism and \(f(e_b)=e_b\),
$$
e_b=f(e_b)\in f(\Ends(U_+))
=\Ends(f(U_+)). 
$$
Applying Fact ~\ref{lem:end-splitting} to \(f(a)\), we know that \(C_b\) is contained either in \(f(U_+)\) or in \(f(U_-)\). \\
If
$
C_b\subset f(U_-),
$
then
$
e_b\in\Ends(C_b)\subset\Ends(f(U_-)),
$ and we get a contradiction.
Hence
$$
C_b\subset f(U_+).
$$
Thus
$$
C_b\subset U_+\quad\Longrightarrow\quad C_b\subset f(U_+).
$$
Exchanging \(+\) and \(-\) gives
$$
C_b\subset U_-\quad\Longrightarrow\quad C_b\subset f(U_-).
$$
Since \(b\) was arbitrary, we obtain
that \(a\) and \(f(a)\) induce the same partition of the boundary components of \(F\).

\end{proof}

\begin{lemma}
Suppose \(a\) is separating. If property (b) holds, then \(a\) and \(f(a)\) induce the same partition of \(\Ends(F)\).
\end{lemma}
\begin{proof}
We now show that, under property~(b), the corresponding partition of the ends of \(F\) is also preserved.
Since \(f\) is a homeomorphism,
$$
\left|\{\text{punctures of }U_+\}\right|
=
\left|\{\text{punctures of }f(U_+)\}\right|.
$$
Denote by 
$
\mathcal B_+$ a collection of boundary components $b$ of 
$\partial F$, such that $C_b\subset U_+$.
Cutting \(U_+\) along the curves in \(\mathcal B_+\), we obtain a decomposition of $U_+$ as:
$$
U_+
=
(F\cap U_+)
\cup
\bigcup_{C_b\subset U_+}C_b.
$$
Hence
$$
\left|\{\text{punctures of }U_+\}\right|
=
\left|\Ends(F)\cap U_+\right|
+
\sum_{C_b\subset U_+}
\left|\{\text{punctures of }C_b\}\right|.
$$
The same decomposition for \(f(U_+)\) gives
$$
\left|\{\text{punctures of }f(U_+)\}\right|
=
\left|\Ends(F)\cap f(U_+)\right|
+
\sum_{C_b\subset U_+}
\left|\{\text{punctures of }C_b\}\right|.
$$
Comparing them, we obtain
$$
\left|\Ends(F)\cap U_+\right|+\Sigma
=
\left|\Ends(F)\cap f(U_+)\right|+\Sigma,
$$
where
$$
\Sigma
=
\sum_{C_b\subset U_+}
\left|\{\text{punctures of }C_b\}\right|.
$$
If \(\Sigma<\infty\), we are done.
Suppose now that \(\Sigma=\infty\). By property~(b), at most one component of \(S\setminus F\) contains punctures. Hence there is a unique component
$
C_{b_0}\subset U_+,
$
containing infinitely many punctures, while all other complementary components contribute only finitely many punctures.
Thus the punctures of \(U_+\) consist of the finite set
$
\Ends(F)\cap U_+
$
together with the countably infinite puncture set contained in \(C_{b_0}\). The same description holds for \(f(U_+)\). Since \(C_{b_0}\) lies on the same side of both \(a\) and \(f(a)\), we can conclude that:
$$
\left|\Ends(F)\cap U_+\right|
=
\left|\Ends(F)\cap f(U_+)\right|.
$$
The same argument applies to \(U_-\). \\
\end{proof}

\begin{lemma}
Suppose \(a\) is separating. The corresponding complementary components of \(F\) have the same genus if properties~(a) and~(b) both hold. 
\end{lemma}
\begin{proof}
One of the two components of \(S\setminus a\) has finite genus. Denote this component by \(U\). Since

$$
f|_U:U\longrightarrow f(U)
$$
is a homeomorphism,
$$
\genus(U)=\genus(f(U)).
$$
Set
$$
W=F\cap U,
\qquad
V=F\cap f(U).
$$
By property~(a), at most one complementary component of \(S\setminus F\) has positive genus. Since \(U\) has finite genus, it cannot contain such a component. Hence every complementary component of \(S\setminus F\) contained in \(U\) has genus zero. Therefore
$$
\genus(U)=\genus(W).
$$
Similarly,
$$
\genus(f(U))=\genus(V).
$$
This ends the proof.\\
\end{proof}
\subsection{The case of at most one accumulated by punctures and at most one non-planar end}
\label{subsec:planar-nonplanar}
Now we focus on the case where both conditions~(P) and~(G) hold. Then we choose $F$ to satisfy properties~(a) and~(b), and Lemmas~3.3,~3.4 and~3.5 above apply. Then we can show:
\begin{lemma}
For every $n$ there exists
$$
[g_n]\in\Map_f(S),
$$
such that
$$
\restr{g_n^{-1}\circ f}{R_n}=id_{R_n}.
$$
\end{lemma}
\begin{proof}
First, we treat a case where $a$ is non-separating. Then $f(a)$ is also non-separating. It follows that the curves belong to the same $\Map(F,\partial F)$-orbit and we get that there exists
$
[g_a]\in\Map(F,\partial F)\subset\Map(S),
$
such that
$
g_a(a)=f(a)
$ and
$
g_a^{-1}\circ f
$
restricts to a homeomorphism of each component of \(S\setminus a\).

If $a$ is separating then combining Lemmas~3.3,~3.4 and~3.5 with Fact~\ref{classification}, we obtain the following: there exists
$
[g_a]\in\Map(F,\partial F)\subset\Map(S),
$
such that
$
g_a(a)=f(a)
$ and
$
g_a^{-1}\circ f
$
restricts to a homeomorphism of each component of \(S\setminus a\).
If \(a_2\) is another boundary component of $P_1$, we apply the same argument to the connected component \(S'\) of \(S\setminus a\) containing \(a_2\), replacing \(f\) by

$$
f'=\left.(g_a^{-1}\circ f)\right|_{S'}.
$$
Repeating this procedure for all boundary components of $P_1$, we obtain
$
[g_1]\in\Map_f(S)
$,
such that
$
g_1(P_1)=f(P_1)
$
and \(g_1\) fixes every boundary component of \(P_1\) set-wise.

A homeomorphism of a pair of pants preserving each boundary component is isotopic, through boundary-preserving homeomorphisms, to one fixing the boundary pointwise. Hence, after modifying \(g_1\) by such an isotopy, we may assume that

$$
f_1:=g_1^{-1}\circ f
$$
restricts to a homeomorphism of \(P_1\) which fixes \(\partial P_1\) pointwise and permutes the punctures of \(P_1\).

Let \(S_2\) be the connected component of \(S\setminus P_1\) containing \(R_2\), and set
$$
f_1'=\left.f_1\right|_{S_2}.
$$
Applying the preceding construction to \(f_1'\), we can choose a finitely supported mapping class $[g_2']\in\Map_f(S_2)$, such that the mapping class of $(g_2')^{-1}\circ f_1'$ fixes \(P_2\).
Then $[f]$ and $[g_2'\circ g_1]$ can only differ by a product $[g_2'']$ of Dehn twists about curves in $\partial{P_1}$ and $\partial{P_2}$. So if we set $g_2=g_2''\circ g_2'\circ g_1$, we get that $[\restr{f}{R_2}]=[\restr{g}{R_2}]\in \Map(R_2, \partial R_2) \subseteq \Map_f(S)$.
We may therefore repeat the construction inductively. \\
\end{proof}
We conclude with the following lemma.
\begin{lemma}
We have
   \begin{equation*}
   \lim_{\longrightarrow}\xi_n=\psi, 
   \end{equation*}
  $ \text{where} \  \psi=[f] \ \text{and} \ \xi_n=[g
_n].$
\end{lemma}
\begin{proof}
By definition of the topology $\xi_k \rightarrow \psi$ if and only if $\forall A$, finite collection of homotopy classes of essential simple closed curves, $\exists k >>0$ s.t. $\xi_k \in \psi U_A$, where $\psi U_A$ is a neighborhood of $\psi$. In our case, if we take any $A$, we can find a finite-type surface $F$, which contains $A$. Moreover, $F$ is contained in $R_n$ for some $n$ and we know that $\forall k>n \ \  \xi_k$ and $\psi$ agree on all curves of $R_n$, so $\xi_k \in \psi U_A$.\\
\end{proof}
\subsection{The case of at most one puncture-accumulated end}
\label{subsec:planar}
We now consider the case when the $(P)$ condition holds, i.e. $S$ has at most one end accumulated by punctures.
Let $a$ be a separating
curve and $F$ a subsurface given by Lemma~3.2; and let $U$ be a component of
$F \setminus a$. If both components of $S\setminus a$ have infinite genus, we cannot guarantee that $\genus(V) = \genus(W)$, where $V = U \cap F$ and $W = f(U) \cap F$.
If this is not the case, there exists a handle shift $[h] \in \Stab(\tilde{E})$, where $h(a)$ induces the same partition of $\partial F$ and $\Ends(F)$ as $a$ and $f(a)$, where $h(a) \subset F$, and such that the genus
of $h(U) \cap F$ and $V$ are equal. Thus, there exists $[g_a] \in \Map_f(F, \partial{F}) < \Stab(\tilde{E})$ such that

\[
g_a \circ h(a) = f(a).
\]

Since $[g_a \circ h] \in \Stab(\tilde{E})$, we can now repeat the argument
above to produce $[k_n] \in \Stab(\tilde{E})$ such that $[k_n] \longrightarrow [f]$.

The next lemma shows that introducing handle shifts is essential in this case.

\begin{lemma}
If $S$ has at least two ends accumulated by genus, then $\Map_f(S)$ is
not dense in $\Stab(\tilde{E})$.
\end{lemma}

\begin{proof}
Choose a separating closed curve $a$ such that each
component of $S \setminus a$ is infinite genus. Let $f$ be an arbitrary element of $\Map_f(S)$.
Let $F$ be an essential finite-type subsurface of $S$ so that
$f \in \Map(F, \partial F) < \Map(S)$ and $F$ satisfies the Lemma~3.2 with the property (b). Let $U$ be a component of $F \setminus a$; then $U$ and $f(U)$
contain the same components of $\partial F$. Further,
\[
\genus(U)=\genus(f(U)),
\]
so either $f([a])=[a]$ or $i([a],f([a]))>0$.

Now, let $b$ be a simple closed curve so that $a$ and $b$ cobound a genus--1
subsurface of $S$. There exists a handle shift
$[h] \in \Stab(\tilde{E})$ with $h(a)=b$. Therefore, no element in $\Map_f(S)$ can agree
with $h$ on $a$;
hence, $\Map_f(S)$ is not dense in $\Stab(\tilde{E})$. 

\end{proof}

\subsection{The case of at most one genus-accumulated end}
\label{subsec:nonplanar}
We now consider the case when $(G)$ condition holds, i.e. $S$ has at most one non-planar end.
Denote by $\Map_p(S)$ the subgroup of $\Stab(\tilde{E})$ generated by puncture shifts and finitely supported mapping classes.
Let $a$ be a separating curve and $F$ is a surface that satisfies Lemma~3.2 with property~(a); and let $U$ be a component of $F \setminus a$. If both components of $S \setminus a$ have infinite number of punctures, we cannot guarantee that $\# \text{punctures}(V)= \# \text{punctures}(W)$, where $W = U \cap F$ and $V = f(U) \cap F$. If this is not the case, there exists a puncture shift $[p]\in \Map_p(S)$, where $p(a)$ induces the same partition of $\partial F$ and of genus as $a$ and $f(a)$, where $p(a) \subset F$, and such that the number of punctures of $p(U) \cap F$ and $V$ is equal. 

Thus, there exists $[g_a] \in \Map(F, \partial F) < \Stab(\tilde{E})$ such that $g_a \circ p(a) = f(a)$ and $[g_a \circ p] \in \Map_p(S)$, so we can now repeat the argument above to produce $[k_n] \in \Map_p(S)$ such that $[k_n] \longrightarrow [f]$.
Moreover, just as before we can show,
\begin{lemma}
If $S$ has at least two ends accumulated by punctures, then $\Map_f(S)$ is
not dense in $\Stab(\tilde{E})$.
\end{lemma}

\subsection{The general case}
We end with the situation where we have no assumptions on $S$. \\
Suppose neither condition $(P)$ or $(G)$ hold, so we cannot guarantee that, if $a$ is separating and $F$ as in Lemma~3.2, $\mathrm{genus}(V) = \mathrm{genus}(W)$
and $\#\mathrm{punctures}(V) = \#\mathrm{punctures}(W)$, where
let $U$ be a component of $F \setminus a$, and set $W := F \cap U$, $V := F \cap f(U)$
as before. But as in the previous sections, we can compose with a handle shift $h$ and/or a puncture shift $p$ to fix the problem. So we can construct $[k_n] \in \langle \mathrm{Map}_f(S), \{\text{handle shifts}\},
  \{\text{puncture shifts}\}\rangle$, such that $[k_n] \to [f]$,
proving density.
Again, the same type of argument shows:
\begin{lemma}
If $S$ has at least two ends accumulated by punctures and at least two non-planar ends, then neither $<\Map_f(S),\{\text{puncture shifts}\}>$, nor $<\Map_f(S),\{\text{handle shifts}\}>$ is
not dense in $\Stab(\tilde{E})$.
\end{lemma}

 This completes the proof of Theorem \ref{maintheorem}.

We finish this section by proving Lemma\ref{prop}. We start with a lemma.
\begin{lemma}
Suppose $S$ has at least one end accumulated by punctures. Let $P$ be any
finite set of punctures of $S$ and $\sigma$ a permutation of $P$. Then
$\sigma$ is realized by a finite composition of puncture shifts, which is finitely supported and fixes every end of $S$ outside
$P$.
\end{lemma}
\begin{proof}
Since every finite permutation is a finite product of transpositions, it suffices to consider a transposition of two punctures $p_0$ and $p_1$. 

Since $S$ has at least one end $\varepsilon$ accumulated by punctures, we may
choose an embedded strip $R\cong\mathbb R\times[-1,1]$ in $S$ such that
\begin{itemize}
  \item $p_0$ and $p_1$ lie in the interior of $R$, and $R$ contains no
        other puncture of $P$;
  \item $R\setminus(\text{punctures of }P)$ contains infinitely many
        punctures, which accumulate to $\varepsilon$ at both sides of $R$.
\end{itemize}
Label the punctures of $R$ along the strip by $q_n$, $n\in\mathbb Z$,
with $q_0=p_0$ and $q_1=p_1$.

Let $h_1$ be the puncture shift along $R$ sending $q_n$ to $q_{n+1}$.
Let $R'$ be a strip obtained from $R$ by pushing its boundary around
$h_1(q_0)=q_1$, so that $R'$ contains all the $q_n$ with $n\neq 1$
and not $q_1$, and let $h_2$ be the puncture shift along $R'$ sending each
$q_n$ ($n\neq 1$) to the next puncture on $R'$, i.e.\ $q_n\mapsto q_{n+1}$
for $n\neq 0,1$ and $q_0\mapsto q_2$. We choose $R'$ and the two shifts
so that $h_2=h_1$ outside of a finite type subsurface  $F\subset S$.

Then $h=h_2^{-1}\circ h_1$ satisfies:
\begin{itemize}
  \item $h(p_0)=p_1$ and $h(p_1)=p_0$, and $h$ fixes every other puncture
        (see Figure~1);
  \item $h$ is the identity outside $F$, so it is supported in a finite-type
        subsurface of $S$; in particular it fixes every end of $S$ other
        than $p_0$ and $p_1$.
\end{itemize}
Thus $h$ realizes the transposition $(p_0\,p_1)$ as a product of two
puncture shifts. Applying this to each factor of a decomposition gives the claim.
\\
\end{proof}
\begin{figure}[ht]
  \centering
  \includegraphics[width=0.8\textwidth]{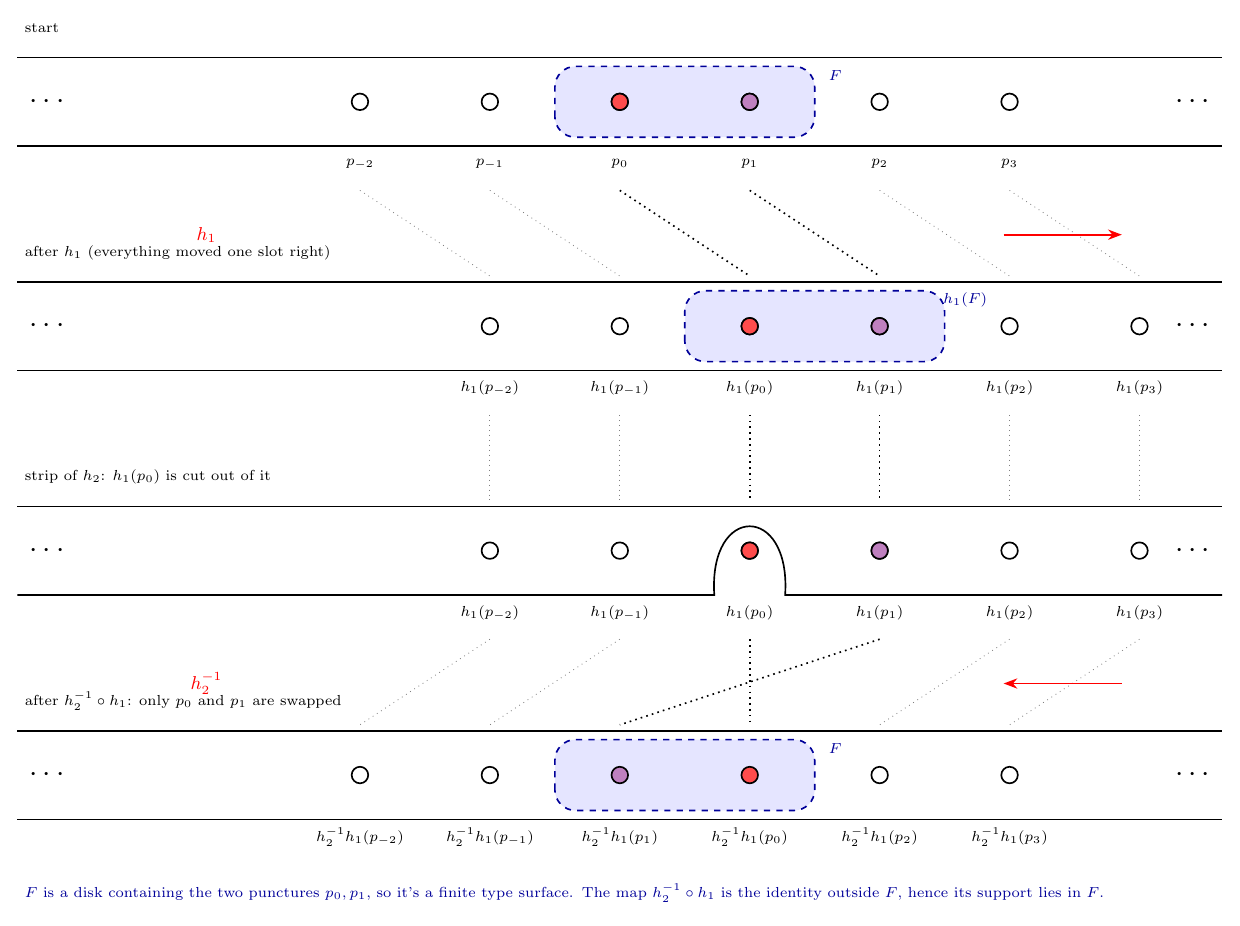}
  \caption{Realizing the transposition of $p_0$ and $p_1$ as $h_2^{-1}\circ h_1$. $K$ is a two-punctured disk.
   The map $h_2^{-1}\circ h_1$ is the identity outside $K$,
   hence $\operatorname{supp}(h_2^{-1}\circ h_1)\subset K$.}
  \label{fig:swap}
\end{figure}
\begin{lemma}
Suppose $S$ has at least one end accumulated by punctures, then
\[
\overline{\operatorname{Map}_{f}(S)}
\subseteq
\overline{
\left\langle
\operatorname{Map}_c(S),\{\text{puncture shifts}\}
\right\rangle}.
\]
\end{lemma}

\begin{proof}
Let $[f]\in\operatorname{Map}_{f}(S)$. Choose a finite-type subsurface $F\subset S$ supporting $f$, and let
\[
P=\{p_1,\ldots,p_n\}
\]
be the punctures contained in $F$. The mapping class $[f]$ induces a
permutation $\sigma_f\in S_n$ of $P$.\\
By the preceding lemma, every finite permutation of punctures
can be realized as a finite composition of puncture shifts. Hence there
exists a composition of puncture shifts $s$ such that
\[
\sigma_s=\sigma_f.
\]
Consequently, the mapping class
\[
[g]=[s^{-1}f]
\]
fixes every puncture in $P$.\\
Thus 
$
[g]\in\operatorname{Map}_c(S).
$
It follows that
\[
[f]=[sg]\in
\left\langle
\operatorname{Map}_c(S),\{\text{puncture shifts}\}
\right\rangle.
\]
Therefore
\[
\operatorname{Map}_{f}(S)
\subseteq
\left\langle
\operatorname{Map}_c(S),\{\text{puncture shifts}\}
\right\rangle,
\]
and taking closures yields
\[
\overline{\operatorname{Map}_{f}(S)}
\subseteq
\overline{
\left\langle
\operatorname{Map}_c(S),\{\text{puncture shifts}\}
\right\rangle}.
\]
\end{proof}

This gives us the proof of Lemma \ref{prop}, together with Corollary \ref{cor1}.

\textbf{Declaration of AI use.} The author used ChatGPT (OpenAI) and Claude (Anthropic) free models as an assistant in preparing this manuscript. It was used for verifying the bibliographic details of cited works, LaTeX assistance, revision of selected mathematical arguments, creation of figures. The authors take full responsibility for all mathematical statements and proofs in this paper. The authors reviewed and verified all changes made with this assistance and take full responsibility for the content of the paper.

\affil[]{\textit{ekaterina.kuritsina@u-pec.fr}}

\affil[]{Univ Paris Est Creteil, Univ Gustave Eiffel, CNRS, LAMA UMR8050, F-94010 Creteil, France}


\bibliographystyle{alpha}
\bibliography{references}
\begin{flushleft}
    Email: \texttt{ekaterina.kuritsina@u-pec.fr}\\
    Univ Paris Est Creteil, Univ Gustave Eiffel, CNRS, LAMA UMR8050, F-94010 Creteil, France
\end{flushleft}

\end{document}